\documentclass[11pt]{amsart}

\usepackage{amsmath,amssymb,amsfonts,enumerate,amsthm, amscd}

\renewcommand{\dim}{\mbox{dim}\,}

\newcommand{\pd}{\mbox{pd}\,}

\newcommand{\Z}{\mathbb{Z}}

\newcommand{\N}{\mathbb{N}}

\newcommand{\fa}{\mathfrak{a}}
\newcommand{\fb}{\mathfrak{b}}
\newcommand{\fm}{\mathfrak{m}}

\newtheorem{theorem}{Theorem}[section]
\newtheorem{lemma}[theorem]{Lemma}
\newtheorem{proposition}[theorem]{Proposition}
\newtheorem{corollary}[theorem]{Corollary}

\theoremstyle{definition}

\theoremstyle{remark}
\newtheorem{remark}[theorem]{Remark}
\newtheorem{example}[theorem]{Example}

\begin{document}

\title[On valuation rings]{On valuation rings}

\author{Mohammed Kabbour}
\address{Mohammed Kabbour\\Department of Mathematics, Faculty of Science and Technology of Fez, Box 2202, University S.M. Ben Abdellah Fez, Morocco.}

\author{Najib Mahdou}
\address{Najib Mahdou\\Department of Mathematics, Faculty of Science and Technology of Fez, Box 2202, University S.M. Ben Abdellah Fez, Morocco.}



\begin{abstract}

In this paper we provide necessary and sufficient conditions for $
R=A\propto E $ to be a valuation ring where $E$ is a non-torsion
or finitely generated $A-$module. Also, we investigate the $ (n,d)
$ property of the valuation ring.

\end{abstract}

\maketitle

\bigskip
\begin{section} {\ Introduction }
\bigskip

All rings considered will be commutative and have identity
element; all
modules will be unital. \\

Let $A$ be a ring, an $A-$module $E$ is said to be
\textit{uniserial} if the set of its submodules is
 totally ordered by inclusion; equivalently, for every $(x,y)\in
E^{2} , x\in Ay$ or $y\in Ax$. A ring $A$ is called a
\textit{valuation ring} if $A$ is an uniserial $A-$module. We note
that $A$ is a valuation ring if and only if $A$ is a local ring
and every finitely generated ideal is principal. See for
instance [\cite{Co1}, \cite{Co2}, \cite{FS}, \cite{H}, \cite{J}, \cite{W}]. \\

An \textit{arithmetical ring} is a ring $A$ for which the ideals
form a distributive lattice, i.e for which $(\fa+\fb)\cap \frak
c=(\fa\cap\frak c)+(\fb\cap\frak c)$
 for all ideals of $A.$ In \cite{J} C.U. Jensen
 gives some more characterization of arithmetical ring, it is proved that a ring $A$ is
 an arithmetical ring if and only if every
 localization $A_\fm$ at a maximal (prime) ideal $\fm$ is a valuation ring. See for instance
 [\cite{BKM}, \cite{BG}, \cite{Co1}, \cite{Co2}, \cite{FS}, \cite{G1}, \cite{J}]. \\

 Let $A$ be a ring, $E$ be an $A$-module and $R:=A\propto E$ be the
set of pairs $(a,e)$ with pairwise addition and multiplication
given by $(a,e)(b,f) = (ab,af+be)$. $R$ is called the trivial ring
extension of $A$ by $E$ (also called the idealization of $E$ over
$A$). Considerable work, part of it summarized in Glaz's book
\cite{G1} and Huckaba's book \cite{H}, has been concerned with
trivial ring extensions. These have proven to be useful in solving
many open problems and conjectures for various contexts in
(commutative and non-commutative) ring theory. See for instance
[\cite{BKM}, \cite{BG}, \cite{G1}, \cite{H}, \cite{KM1}, \cite{KM2}, \cite{PR}]. \\

 For nonnegative integer $n,$ an $A-$module $E$ is said to be of finite
$n-$presentation ( or $n-$presented ) if there exists an exact
sequence:$$F_n \rightarrow F_{n-1} \rightarrow ...\rightarrow F_1
\rightarrow F_0 \rightarrow E \rightarrow 0$$ where $F_i$ is a
free $A-$module of finite rank. We write $$\lambda_A (E)=\sup
\left \{ n;\mbox {there exists a finite n-presentation of } E
\right\}.$$
 The $\lambda-$dimension of a ring $A\left (\lambda-\dim A\right)$
 is the least integer $n$ ( or $\infty$ if none such exists )
 such that $\lambda_A (E)\geq n$ implies $\lambda_A (E)=\infty $. $A$ is called a strong
 $n$-coherent ring in [\cite{KM1}, \cite{KM2}, \cite{M1}], \cite{M2}]. \\
Throughout, $\pd_A (E)$ will denote the projective dimension of
$E$ as an $A-$module. \\

Given nonnegative integers $n$ and $d$, we say that a ring $R$ is
an $(n,d)$-ring if $pd_{R}(E) \leq d$ for each $n$-presented
$R$-module $E$ (as usual, $pd$ denotes projective dimension). For
integers $n$, $d \geq  0$ Costa asks in \cite{C} whether there is
an $(n,d)$-ring which is neither an $(n,d-1)$-ring nor an
$(n-1,d)$-ring? The answer is affirmative for $(0,d)$-rings,
$(1,d)$-rings, $(2,d)$-rings and $(3,d)$-rings for each integer
$d$. See for instance [\cite{C}, \cite{CK}, \cite{KM1},
\cite{KM2}, \cite{M1}], \cite{M2}, \cite{Z}]. \\

The goal of section 2 of this paper is to provide necessary and
sufficient conditions for $ R :=A\propto E $ to be a valuation
ring where $E$ is a non-torsion or finitely generated $A-$module.
The section 3 is devoted to investigate the $ (n,d)$-property of
the valuation ring.
\end{section}

\bigskip
\begin{section}{\ Trivial extensions defined by valuation ring}
\bigskip

This section develops a result of the transfer of valuation
property to trivial ring extension.  Recall that an $A$-module $E$
is called a torsion module if for every $u \in E$, there exists $0 \not= a \in A$
such that $au =0$.  \\

\bigskip

\begin{theorem}\label{pr}

Let $A$ be a ring and $E$ an nonzero $A-$module. Let $R :=A\propto
E$ be the trivial ring extension of $A$ by $E.$
\begin{enumerate}
    \item Assume that $E$ is a non-torsion $A-$module. Then $R$ is a
    valuation ring if and only $A$ is a valuation domain and $E$ is
    isomorphic to $K :=qf(A),$ the field of fractions of $A.$
    \item Assume that $E$ is a finitely generated $A-$module. Then $R$ is
    a valuation ring if and only if $A$ a is field and $E\simeq A.$
\end{enumerate}
\end{theorem}

\bigskip

Before proving { Theorem } \ref {pr}, we establish the following
Lemma:

\bigskip

\begin{lemma}\label{vr}  Let $A$ be a
ring, $E$ a non zero $A-$module and let $R: =A\propto E$ be the
trivial ring extension of $A$ by $E.$ If $R$ is a valuation ring
then $A$ is a valuation domain and $E$ is an uniserial $A-$module.
\end{lemma}

\begin{proof}  Assume that $R$
is a valuation ring. First we wish to show that $A$ is a valuation
ring and $E$ is a uniserial $A-$module. Let $(a,b)\in A^{2}$, if
$(a,0)$ divides $(b,0)$ (resp., $(b,0)$ divides $(a,0)$) then $a$
divides $b$ (resp., $b$ divides $a$). Hence $A$ is a valuation
ring. On the other hand, let $(x,y)\in E^{2}$. If $(0,x)$ divides
$(0,y)$ (resp., $(0,y)$ divides $(0,x)$) then there exists $(c,z)
\in R$ such that $(0,y) =(c,z)(0,x)$ (resp., $(0,x) =(c,z)(0,y)$)
and so $y \in Ax$ (resp., $x \in Ay$). Therefore $E$ is an
uniserial $A-$module.\\

 We claim that $A$ is an integral domain. Deny. Let$(a,b)\in A^{2}$ such
that $ab=0,a\neq 0 $ and $ b\not=0.$ For each $x\in E$, $(b,0)$
divides $(0,x)$ (since $R$ is a valuation ring and $(0,x)$ does
not divides $(b,0)$ (since $b \not= 0$)) and so there exists $y\in
E$ such that $by=x$, thus $ax=0$ and so $a\in (0:E)$. Also, for
each $x \in E$, $(a,0)$ divides $(0,x)$ and so $x \in aE =0$, a
contradiction since $E \not= 0$. Thus $A$ is an
integral domain.\\

\end{proof}

\bigskip

\begin{proof} \textit{of Theorem 2.1.} \\
 {\bf 1)} Assume
that $A$ is a valuation domain and let $R :=A \propto K$, where $K
:=qf(A)$. Our aim is to show that $R$ is a valuation ring. Let
$(a,x)$, $(b,y) \in R-\{0,0\}$. Two cases are then possibles: \\
{\bf Case 1.} $a =b =0$. There exists then $c \in A$ such that $x
=cy$ (resp., $y =cx$) since $K :=qf(A)$ and $A$ is a valuation
domain. Hence, $(0,x) =(c,0)(0,y)$ (resp., $(0,y) =(c,0)(0,x)$) as
desired. \\
{\bf Case 2.} $a \not= 0$ or $b \not= 0$. We may assume that $a
\not= 0$ and $b \in Aa$. Let $c \in A$ such that $ac =b$ and let
$z \in K$ such that $az + cx =y$. Hence, $(a,x)(c,z) =(b,y)$ as
desired. \\

Conversely, assume that $E$ is a non-torsion $A$-module and $R$ is
a valuation ring. We wish to show that $E\simeq K.$ Let $u\in E$
such that $ ( 0:u)=0 $ and let $f:K\otimes Au \rightarrow K\otimes
E$ be the homomorphism of $ A- $module induced by the inclusion
map $ Au\hookrightarrow E. $ Since the field $ K $ is a flat
$A-$module, then $ f $ is injective. Let $ (\lambda,x)\in K\times
E, $ by {Lemma} \ref{vr} we get that $ x\in Au $ or $ u\in Ax. $
If $ x=au $ for some $ a\in A $ then $f(\lambda\otimes
au)=\lambda\otimes x. $ If $ u\in Ax $ then there exists $ a\in A
$ such that $ u=ax. $ Thus
$$f\left(\frac{\lambda}{a}\otimes
u\right)=\frac{\lambda}{a}\otimes u=\frac{\lambda}{a}\otimes
ax=\lambda\otimes x.$$ Consequently, $f$ is an isomorphism of
$A-$module. Now,  consider the homomorphism of $ A-$module $
g:E\rightarrow K\otimes E $ defined by $ g(x)=1\otimes x$. For all
multiplicatively closed subset $ S $ of $ A, $ the $
S^{-1}A-$modules $ S^{-1}E $ and $ S^{-1}A\otimes_{A}E $ are
isomorphic; more precisely the map $ \varphi : S^{-1}E\rightarrow
S^{-1}A\otimes_{A}E,$ where $\varphi\left(\displaystyle
\frac{x}{s}\right)=\displaystyle\frac{1}{s}\otimes x$ is
isomorphism. If $g(x)=1\otimes x=0,$ then there exists $0\neq a\in
A$ such that $ax=0.$ By Lemma \ref{vr} $ x\in Au $ or $ u\in Ax$.
But $u \notin Ax$ since $ax =0$, $a \not= 0$ and $(0:u) =0$.
Hence, $ x=bu $ for some $ b\in A $.  Then $ abu=0, $ hence $ ab=0
$ since $(0:u) =0$ and so $ b=0$ (since $A$ is a valuation domain
and $a \not= 0$); thus $ x=0. $ It follows that $g$ is injectif.
Let $(\lambda,x)\in K\times E,$ if $\lambda \in A$ then
$\lambda\otimes x=1 \otimes\lambda x=g(\lambda x).$ Now if $
\lambda^{-1}\in A $ then there exists $ y\in E $ such that $
\lambda^{-1}y=x, $ since $ (\lambda^{-1},0) $ divides $ (0,x). $
Hence
$$\lambda\otimes x=\lambda \otimes(\lambda^{-1}y)=1\otimes y=g(y).$$
Consequently, $g$ is an isomorphism of $A-$module. We deduce that
$$ E\simeq K\otimes_{A} E \simeq K\otimes_{A} Au \simeq K\otimes_{A} A \simeq K.$$

{\bf 2}) If $A$ is a field, then $R :=A \propto A$ is a valuation
ring by the proof of {\bf 1)} above. Conversely, assume that $ E $
is a finitely generated $A-$module. We denote by $ \fm $ the
maximal ideal of $ A. $ By Lemma \ref{vr}, $ E/\fm E $ is an $
A/\fm-$vector space and for all $ (x,y)\in
E^{2},\bar{x}\in(A/\fm)\bar{y} $ or $ \bar{y}\in(A/\fm)\bar{x}. $
Hence, $ \dim_{A/\fm A}\left(E/\fm E\right) =0 $ or 1. If $ E=\fm
E, $ then $ E=0 $ by Nakayama Lemma which is absurd. Thus $ E/\fm
E=\left(A/\fm \right)\bar{v} $ for some $ v\in E\setminus \fm E. $
By Nakayama Lemma $ v $ generate $ E. $ Suppose that $ \fm \neq0,
$ let $ 0\neq a\in\fm, $ we have $ (a,0) $ divides $ (0,v) $, then
there exists $ b\in A $ such that $ (a,0)(0,bv)=(0,v). $ Hence $
(1-ab)v=0, $ therefore $ v=0$ (since $1-ab$ is unit), which is
absurd. Therefore, the ring $A$ is a field and $ E=Av, $ therefore
$ E\simeq A$ completing the proof of Theorem \ref{pr}. \\
\end{proof}

\bigskip

\begin{corollary}
Let $ A $ be a ring, then $ A\propto A $ is a valuation ring if
and only if $ A $ is a field.
\end{corollary}

\bigskip

 Theorem \ref{pr} enriches the literature with new examples
of valuation rings, as shown below.

\bigskip

\begin{example}
 Let $ k $ be a field. Let $ A=k[[x]] $ the ring of formal
 power series with coefficients in $ k $ and $ K=k((x)) $ its field
 of fractions. The trivial ring extension of $ A $ by $ K, A\propto
 K $ is a valuation ring.
\end{example}

\bigskip

\begin{example}
    Let $\mathbb{Q}_{p}$ be the completion of $\mathbb{Q}$
    in the $p-$adic topology where $p$ is prime integer. The ring
    of $p-$adic integers is
    $\Z_{p}=\{x\in\mathbb{Q}_{p}:|x|_{p}\leq1\},\mathbb{Q}_{p}$
    is its field of fractions. Then the trivial ring extension of
    $\Z_{p}$ by $\mathbb{Q}_{p}$ is a valuation ring.
\end{example}

\bigskip

Now, we are able to construct a non-valuation ring.

\bigskip

\begin{corollary}
     Let $A$ be a ring and let $E$ be a mixed module, i.e $E$
    is neither torsion nor torsion-free. Then $A\propto E$ is not
    a valuation ring.
\end{corollary}

\begin{proof} Assume that $ A\propto E $ is
a valuation ring and that $ E $ is a non torsion $ A-$module. Let
$ u\in E $ such that $ (0:u)=0. $ By Lemma \ref{vr}, for each $
0\neq x\in E,x\in Au \mbox{ or } u\in Ax$. Suppose that $ u=ax $
for some $ a\in A$. Since $(0:u) =0$, the following implications
hold:
$$ \alpha x=0\Rightarrow\alpha
    ax=0\Rightarrow\alpha
    u=0\Rightarrow\alpha=0. $$ Hence $ (0:x)=0. $ Finally, it is easy to get
    the equality $ (0:x)=0 $ in the case $x \in Au$ since $A$ is an integral
    domain (by Lemma \ref{vr}). Thus $ E $ is
    a torsion-free $A-$ module.
\end{proof}

\end {section}

\bigskip
\begin {section}{\ $(n,d)$-properties and valuation rings}
\bigskip

 Let $A$ be a ring. An $A-$ module is called a cyclically
presented module if it is isomorphic to $A/aA$ for some $a\in A.$
\\
Now, we are able to give our main result in this section. \\

\begin{theorem}\label{ndvr}
Let $ A $ be a valuation ring and let $ Z $ be the subset of its
zero divisors.
\begin{enumerate}
    \item If $ (0:a) $ is not a finitely generated ideal for every
    $ a\in Z\setminus 0 $, then $ A $ is a (2,1)-ring.
    \item If $ (0:a) $ is a finitely generated ideal for some
    $ a\in Z\setminus 0 $, then $ A $ is not a (2,d)-ring for every nonnegative
    integer $ d. $
\end{enumerate}
\end {theorem}

\bigskip

In ordre to prove Theorem \ref {ndvr}, we will use the following
Lemma:

\bigskip

\begin {lemma} ( \cite[Theorem 1]{W}\label{warf}) \\A finitely presented module over a
valuation ring is a finite direct sum of cyclically presented
modules.
\end {lemma}

\bigskip

\begin{proof} \textit{of Theorem 3.1.} \\
{\bf 1)} Assume that $ (0:a) $ is not a finitely generated ideal
for every $ a\in Z\setminus 0 $ and let $ E $ be a 2-presented $
A-$module. We want to show
    that $ pd_{A}(E)\leq1. $ By the above Lemma \ref {warf}, $ E $ is finite direct
    sum of cyclically presented modules; i.e $ E=\bigoplus \limits_{i=1}^{n}Ax_{i}$ and $Ax_{i}\simeq A/a_{i}A $ for some $ a_{i}\in A. $
    Consider the following
    exact sequences: $$ 0\rightarrow a_i A \rightarrow A \rightarrow A/a_{i}A \rightarrow 0
    \eqno{(1)} $$  $$ 0\rightarrow \left(0:a_i\right) \rightarrow A \rightarrow a_{i}A \rightarrow
    0.
    \eqno{(2)} $$ Then $ a_i A $ is a finitely presented $ A-$module and
    $ \left(0:a_i\right) $ is a finitely generated ideal (since $ A/a_{i}A $ is a 2-presented $ A-$module). Therefore, $ a_{i}=0 $ or
    $ a_{i} $ is a nonzero divisor element of $ A$ by hypothesis. Using the exact
    sequence (1), we can also deduce that $ pd\left( A/a_i A \right) \leq 1$ and so $ pd(E)=\sup\left\{pd_A
\left( A/a_i A\right):1\leq i\leq n
    \right\} \leq 1$. Hence, $ A $ is a $ (2,1)-$ring. \\

    {\bf 2)} Let $ a\in Z\setminus 0 $ such that $ (0:a) $ is a finitely
    generated ideal. The following exact sequences of $ A-$modules:
    $$ 0\rightarrow aA \rightarrow A \rightarrow A/aA \rightarrow 0
    \eqno{(1)} $$
     $$ 0\rightarrow (0:a) \rightarrow A \rightarrow aA \rightarrow 0
    \eqno{(2)} $$ show that $ aA $ is a 1-presented $ A-$module and
    $ A/aA $ is a 2-presented $A-$module. But the $ \lambda -$dimension
    of every valuation ring is at most two (\cite[Corollary 2.12]{Co1}), hence $ \lambda_A\left(A/aA\right)=\infty . $
    Let $ b\in A $ such that $ (0:a)=bA. $ By the following exact
    sequence of $ A-$modules:
  $$ 0\rightarrow (0:b) \rightarrow A \rightarrow bA \rightarrow 0
    \eqno{(3)} $$ we get that $ (0:b) $ is finitely generated. Then
    there exists $ c\in A $ such that $ (0:b)=cA. $ We claim that
    $ (0:c)=bA$. \\
     Indeed, $ bA\subseteq (0:c) $  since $ bc=0$. On the other hand, let
    $ x\in (0:c) $, that is $ cx=0$. But $a \in (0:b) =cA$ (since
    $(0:a) =bA$), then $a =ct$ for some $t \in A$. Hence, $ax =cxt
    =0$ and so $x \in (0:a) =bA$, as desired. \\
    By using the exact sequences (3) and $$ 0\rightarrow (0:c) :=bA  \rightarrow A \rightarrow cA \rightarrow 0
    \eqno{(4)} $$ we get the equalities
    $ pd_A\left(bA\right)=pd_A\left(cA\right)+1 $ and
    $ pd_A\left(cA\right)=pd_A\left(bA\right)+1. $ Hence
    $ pd_A\left(bA\right)=pd_A\left(cA\right)=\infty. $ Therefore
    $ pd_A\left(A/aA\right)=\infty, $ which completes the proof
    of Theorem \ref{ndvr}.
\end{proof}

\bigskip

Now, we construct a $(2,1)$-ring which is a particular case of
\cite[Theorem 3.1]{KM2}. \\

\bigskip

\begin{corollary}
Let $ A $ be a valuation domain which is not a field, $ K :=qf(A)$
and let $ R :=A\propto K $ be the trivial ring extension of $ A $
by $ K. $ Then $ R $ is a $ (2,1)-$ring.
\end{corollary}

\begin{proof}
By {Theorem} \ref{pr} $ R $ is a valuation ring. Let $ 0\not=
(a,x)\in R. $ It is easy to get successively that $ (a,x) $ is a
zero divisor if and only if $ a=0, $ and the equality $
(0:(0,x))=0\propto K. $ Assume that there exists $0 \not= x \in K$
such that $ (0:(0,x)) $ is not a finitely generated ideal. Then
there exists $ x_1,...,x_n\in K $ such that
$$ (0:(0,x))=R(0,x_1)+...+R(0,x_n)=0\propto
\left(Ax_1+...+Ax_n\right). $$ Therefore, $ K=Ax_1+...+Ax_n. $ We
put $ x_i=\displaystyle \frac{a_i}{d} $ for every $ 1\leq i\leq n,
$ where $ a_i\in A $ and $ 0\not= d\in A. $ Hence, $K =dK =Aa_{1}
+ \ldots + Aa_{n} \subseteq A$, a contradiction. Therefore, $ R $
is a $ (2,1)-$ring.
\end{proof}

\bigskip

Let $ A $ be a ring, a necessary and sufficient conditions for $ A
$ to be coherent is that $ (0:a) $ is finitely generated ideal for
every element $ a\in A, $ and the intersection of two finitely
generated ideals of $ A $ is a finitely generated ideal of $ A. $
By {Theorem} \ref{ndvr}, we have: \\

\bigskip

\begin{corollary}
Let $ A $ be a valuation ring with zero divisors, if $ A $ is
coherent then $ A $ is not a $ (2,d)-$ring for every nonnegative
integer $ d. $
\end{corollary}

\bigskip

\begin{remark}
Let $A$ be a valuation ring and $\fm$ its maximal ideal. By
\cite[proposition 2.10]{Co1}, $A$ is $(2,1)-$ring if and only if
$\fm$ is flat.
\end{remark}

\bigskip

\begin{example}
Let $p$ be a prime nonnegative integer and $n\in  \N^*$. The
valuation ring $ \Z /p^n \Z $ is not a $ (2,d)-$ring for every
nonnegative integer $ d. $
\end{example}

\bigskip

\begin{example}
Let $ K $ be a field and $ n $ an integer such that $ n\geq 2. $
We denote $ A=K[x]/\left(x^n\right) $ and $ \overline{P}=P+A $ for
every $ P\in K[x]. $ It is easy to see that $ A $ is a valuation
ring. We have $ \left(0:\overline{x^{n-1}}\right)=\overline{x}A, $
which is a finitely generated ideal. Hence $ A $ is not a $
(2,d)-$ring for every positive integer $ d. $
\end{example}

\bigskip

Now, we study the relationship between the $(n,d)$-properties and
an arithmetical rings. \\

\bigskip

\begin{proposition}
Let $A$ be an arithmetical ring.\begin{enumerate}
    \item Suppose that for every maximal ideal $ \fm $ of $ A $ and $ 0\not=x\in
    Z\left(A_\fm\right),\,(0:x) $ is not a finitely generated
    ideal. Then $A$ is a $ (3,1)-$ring.
     \item If $ (0:x) $ is finitely generated for some maximal ideal
     $ \fm $ of $ A $ and $ 0\not=x\in
    Z\left(A_\fm\right) $, then $ A $ is not a $(1,d)-$ring for every
    nonnegative integer $ d. $
\end{enumerate}
\end{proposition}

\begin{proof}
(1) Let $ E $ be a 3-presented $ A-$module   and let $ \fm $ be a
maximal ideal of $ A. $ The $ A_\fm -$module $ E_\fm $ is then
2-presented. For every $ x\in Z\left(A_\fm \right),\,(0:x) $ is
not finitely generated and $ A_\fm $ is a valuation ring, then by
Theorem \ref{ndvr} $ A_\fm $ is $ (2,1)-$ring. Thus the projective
dimension of $ E_\fm $ over $ A_\fm $ is at most one. Since $
\lambda-\dim A\leq3 $ ( by \cite[Theorem 2.1]{Co1}), then $ E $
admits a finite free resolution. We deduce that
$$ pd_A(E)=\sup\left\{pd_{A_\fm}\left(E_\fm\right);\,\fm\, \mbox
{ is a maximal ideal of } A\right\} \leq 1$$ and so $ A $ is a $
(3,1)$-ring.
\par
(2) Assume that $ (0:x) $ is a finitely generated ideal for some
maximal ideal $ \fm $ of $ A $ and $0 \not= x\in
Z\left(A_\fm\right). $ We put $ x=\displaystyle\frac{a}{s} $ where
$ a\in A $ and $ s\notin \fm. $   Then $
pd_{A_\fm}\left(A_\fm/aA_\fm\right)=pd_{A_\fm}\left((A/aA)_\fm\right)$
since the $ A_\fm-$modules $ A_\fm/aA_\fm $ and $
\left(A/aA\right)_\fm $ are isomorphic. Hence, $
pd_{A_\fm}\left((A/aA)_\fm\right)=\infty $ by {Theorem}
\ref{ndvr}(2)and so $ pd_A(M)\geq pd_{A_\fm}\left(M_\fm\right) $
for every $ A-$module $ M, $ then $
pd_{A}\left(A/aA\right)=\infty. $ On the other hand, $ A/aA $ is a
finitely presented $ A-$module. Consequently, $ A $ is not a $
(1,d)$-ring for every nonnegative integer $ d. $
\end{proof}

\end {section}
\bigskip




\bigskip\bigskip


\end {document}